\documentclass[11pt]{amsart}
\usepackage{amsmath}
\usepackage{amsfonts}
\usepackage{latexsym}
\usepackage{amssymb}
\usepackage{enumerate}
\usepackage[dvips]{graphics}

\newcommand{\Z}{{\mathbf Z}}

\newcommand{\R}{{\mathbf R}}
\newcommand{\C}{{\mathbf C}}

\newcommand{\kk}{{\mathbf k}}
\newcommand{\cat}{{\rm {cat }}}
\newcommand{\Cat}{\rm {Cat}}

\newcommand{\Hom}{{\rm Hom}}
\newcommand{\im}{\rm im}

\newcommand{\Int}{{\rm {Int}}}

\newcommand{\cwgt}{{\rm {cwgt}}}
\newcommand{\ccat}{{\rm {ccat}}}

\newcommand{\RP}{\mathbf {RP}}
\newcommand{\swgt}{{\rm {swgt}}}

\newcommand{\rk}{{\rm rk}}
\newcommand{\comment}[1]{}

\def\V{\mathcal V}

\def\ker{{\rm Ker }}
\def\im{{\rm Im }}

\def\R{\mathbf R}

\newtheorem{theorem}{Theorem}
\newtheorem{proposition}{Proposition}
\newtheorem{lemma}[proposition]{Lemma}
\newtheorem{corollary}[proposition]{Corollary}
\theoremstyle{definition}
\newtheorem{definition}{Definition}
\newtheorem{example}{Example}
\newtheorem{remark}{Remark}

\begin{document}
\title[Category weights and estimates for {\large {$\cat^1(X,\xi)$}}]{Homological category weights \\
and
estimates for  {\Large $\cat^1(X,\xi)$}}

\author[M. ~Farber and D. ~Sch\"utz]{M.~Farber and D. ~Sch\"utz}
\address{Department of Mathematics, University of Durham, Durham DH1 3LE, UK}
\email{Michael.Farber@durham.ac.uk}

\address{Department of Mathematics, University of Durham, Durham DH1 3LE, UK}
 \email{dirk.schuetz@durham.ac.uk}



\date{}


\begin{abstract} In this paper we study a new notion of category
weight of homology classes developing further the ideas of
 E. Fadell and S. Husseini \cite{FH}.
In the case of closed smooth manifolds the homological category
weight is equivalent to the cohomological category weight of E.
Fadell and S. Husseini but these two notions are distinct already
for Poincar\'e complexes. An important advantage of the
homological category weight is its homotopy invariance. We use the
notion of homological category weight to study various
generalizations of the Lusternik - Schnirelmann category which
appeared in the theory of closed one-forms and have applications
in dynamics. Our primary goal is to compare two such invariants
$\cat(X,\xi)$ and $\cat^1(X,\xi)$ which are defined similarly with
reversion of the order of quantifiers. We compute these invariants
explicitly for products of surfaces and show that they may differ
by an arbitrarily large quantity. The proof of one of our main
results, Theorem \ref{main2}, uses an algebraic characterization
of homology classes $z\in H_i(\tilde X;\Z)$ (where $\tilde X\to X$
is a free abelian covering) which are movable to infinity of
$\tilde X$ with respect to a prescribed cohomology class $\xi\in
H^1(X;\R)$. This result is established in Part II which can be
read independently of the rest of the paper.
\end{abstract}

\maketitle

\section{Introduction}

In this paper we study various generalizations of the classical
Lusternik - Schnirelmann category $\cat(X)$ which arise in
topology of closed one-forms. They are homotopy invariants of
pairs $(X, \xi)$ where $X$ is a finite polyhedron and $\xi\in
H^1(X;\R)$ is a real cohomology class. Several potentially
different notions
\begin{eqnarray}\label{all} \cat(X,\xi) \leq \cat^1(X,\xi)
\leq \Cat(X,\xi)\end{eqnarray} play different roles in application
of the theory of closed one-forms to dynamics, see \cite{farbe2},
\cite{FK}, \cite{farbe3}; each of these invariants turns into the
classical $\cat(X)$ when $\xi=0$. One of the objectives of the
present paper is to show that $\cat^1(X,\xi)$ can be distinct from
$\cat(X,\xi)$ and moreover their difference can be arbitrarily
large. At the moment we have no examples where $\Cat(X,\xi)$ is
distinct from $\cat^1(X,\xi)$.

It is well-known that a most effective lower bound for the
classical Lusternik - Schni\-rel\-mann category $\cat(X)$ is the
cohomological cup-length, i.e. the largest number of cohomology
classes of positive degree such that their cup-product is
nontrivial. In our recent preprint \cite{FS1} we established
cohomological cup-length type lower bounds for $\cat(X,\xi)$ which
use local systems of a special kind. In view of (\ref{all}) all
lower bounds for $\cat(X,\xi)$ hold for $ \cat^1(X,\xi)$ as well.
In order to distinguish between these two invariants one needs to
have lower bounds for $\cat^1(X,\xi)$ which in general are not
true for $\cat(X,\xi)$. Such lower bounds are found in the present
paper.

Our main results are based on the idea of category weight which
was initially introduced by E. Fadell and S. Husseini who proposed
in \cite{FH} to attach \lq\lq weights\rq\rq\, to cohomology
classes so that classes of higher weight contribute more into the
cup-length estimate; see \S \ref{basdef} for more detail. We would
like to mention also papers of Y. Rudyak \cite{Ru} and J. Strom
\cite{Str} who suggested a useful modification of this notion. In
this paper we propose yet another variation of this idea: we
attach weights to homology classes (and not to cohomology classes
as did the previous authors) and measure \lq\lq the level of
nonvanishing\rq\rq\, of a cup-product $u_1\cup u_2\cup \dots\cup
u_r$ 
by evaluating it on homology classes $\langle u_1\cup u_2\cup
\dots\cup u_r, z\rangle$ of different weight. We show that the
notion of category weight of homology classes has an important
advantage of being homotopy invariant (unlike the weights of
Fadell and Husseini). We prove that for closed manifolds the
category weight of a homology class equals the category weight of
Fadell and Husseini of the dual cohomology class. We also show
that this statement is false for Poincar\'e complexes. The results
about category weights of homology classes occupy Part I which can
be read independently of the rest of the paper.

Part II also covers a story which may be read independently of
Parts I and III. Here we study free abelian covers $p: \tilde X
\to X$ and homology classes $z\in H_i(\tilde X;\Z)$ which can be
realized by singular cycles lying arbitrarily far in a specified
direction. Such \lq\lq directions\rq\rq\, are parametrized by
cohomology classes $\xi\in H^1(X;\R)$ with $p^\ast(\xi)=0$. Our
result states that this property of $z$ is equivalent to the
existence of an infinite chain $c'$ such that $\partial c'=c$ and
$c'$ is {\it \lq\lq automatically produced out of finite
data\rq\rq,}\, see the discussion after Theorem \ref{movability}.
The main result of Part II generalizes Theorem 5.3 of
\cite{farbe4} which treats the case of rank one cohomology
classes. It also generalizes our previous result \cite{FS}
covering the case of homology classes with coefficients in a
field; in \cite{FS} our arguments use a different algebraic
mechanism which fails to work over the integers.

In Part III we use the results of Parts I and II to obtain new
cohomological lower bounds for $\cat^1(X, \xi)$. Our Theorem
\ref{main2} gives in many cases stronger estimates than Theorem
5.6 of \cite{farbe4}; note that the latter theorem applies only in
the special case of rank one cohomology classes although the
results of the present paper are valid in full generality and do
not impose this restriction. In Part III we also introduce  a
controlled version of $\cat^1(X,\xi)$ which behaves better under
cartesian products. Finally we compute $\cat^1(X,\xi)$ for
products of surfaces as function of the cohomology class $\xi\in
H^1(X;\R)$. We compare our results with the computations of the
invariant $\cat(X, \xi)$ completed in \cite{FS1}. We conclude that
$\cat^1(X,\xi)$ may exceed $\cat(X,\xi)$ by an arbitrary large
amount.

The following diagram illustrates dependence of parts of this
paper:
$$\begin{array}{ccccc}
\mbox{Part I} & &&& \mbox{Part II}\\
& \searrow & &\swarrow &\\
&& \mbox{Part III} &&
\end{array}
$$
Parts I and II can be read independently, the results of Parts I
\& II are used in Part III.

\specialsection*{\bf Part I: Category weights of homology classes}

Here we introduce and study the notion of category weight of
homology classes which is somewhat dual to the cohomological
notion introduced by E. Fadell and S. Husseini \cite{FH}; the
homological category weight has the advantage of being homotopy
invariant. In part III we use this notion to obtain improved
cohomological lower bounds for $\cat^1(X;\xi)$.

\section{Basic definitions}\label{basdef}

The classical cohomological lower bound for the Lusternik -
Schnirelmann category $\cat(X)$ states that $\cat(X)>n$  if there
exist $n$ cohomology classes of positive degree $u_j\in
H^\ast(X;R_j)$ where $j=1, 2, \dots, n$, such that their
cup-product $u_1 u_2 \dots u_n\not=0 \in H^\ast(X; R)$ is
nontrivial. Here $R_j$ denotes a local coefficient system on $X$
and $R$ is the tensor product $R=R_1\otimes \dots \otimes R_n$.

E. Fadell and S. Husseini \cite{FH} improved this estimate by introducing the
notion of a {\it category weight} $\cwgt(u)$ of a cohomology class $u\in
H^q(X;R)$. Here is their definition:

\begin{definition}\label{def1}
Let $u\in H^q(X; R)$ be a nonzero cohomology class where $R$ is a
local coefficient system on $X$. One says that $\cwgt(u)\geq k$
(where $k\geq 0$ is an integer) if for any closed subset $A\subset
X$ with $\cat_XA\leq k$ one has $u|_A=0\in H^q(A;R)$.
\end{definition}

Recall that the inequality $\cat_XA\leq k$ means that $A$ can be
covered by $k$ open subsets $U_i\subset X$ such that each
inclusion $U_i\subset X$ is null-homotopic, $i=1, \dots, k$.

According to Definition \ref{def1} one has $\cwgt(u) \geq 0$ in
general and $\cwgt(u)\geq 1$ for any nonzero cohomology class of
positive degree. As Fadell and Husseini \cite{FH} showed,
$\cwgt(u)>1$ in some special situations which allows to improve
the lower estimate for $\cat(X)$. Indeed, one has
$$\cat(X) \geq 1
+\sum_{i=1}^n \cwgt(u_i)$$ assuming that the cup-product $u_1u_2\dots
u_n\not=0$ is nonzero.

Y. Rudyak \cite{Ru} and J. Strom \cite{Str} studied a modification
of $\cwgt(u)$, called the strict category weight ${\rm
{swgt}}(u)$. The latter has advantage of being homotopy invariant.
However in some examples the strict category weight is
considerably smaller than the original category weight of Fadell
and Husseini.

In this paper we introduce and exploit a \lq\lq dual\rq\rq\,
notion of category weight of homology classes. It has the
geometric simplicity and clarity of category weight as defined by
Fadell and Husseini but has a surprising advantage of being
homotopy invariant.

\begin{definition}
Let $z\in H_q(X;R)$ be a singular homology class with coefficients in a local
system $R$ and let $k\geq 0$ be a nonnegative integer. We say that
$\cwgt(z)\geq k$ if for any closed subset $A\subset X$ with $\cat_XA\leq k$
there exists a singular cycle $c$ in $X-A$ representing $z$. We say that
$\cwgt(z)=k$ iff $\cwgt(z)\geq k$ and $\cwgt(z)\not\geq k+1$
\end{definition}

In other words, ${\cwgt}(z)\geq k$ is equivalent to the fact that $z$ can be
realized by a singular cycle avoiding any prescribed closed subset $A\subset X$
with $\cat_XA\leq k$.

For example, $\cwgt(z)\geq 1$ iff $z$ can be realized by a singular cycle
avoiding any closed subset $A\subset X$ such that the inclusion $A\to X$ is
homotopic to a constant map.

It will be convenient to define the category weight of the zero homology class
as $+\infty$.

Formally $\cwgt(z)\geq k$ if $z$ lies in the intersection
$$\bigcap_A \im[H_q(X-A;R)\to H_q(X;R)]$$
where $A\subset X$ runs over all closed subsets with $\cat_XA\leq k$.

The relation $\cwgt(z) \leq k$ means that there exists a closed subset
$A\subset X$ with $\cat_XA \leq k+1$ such that any geometric realization of $z$
intersects $A$. In particular we obtain the following inequality
\begin{eqnarray}\label{one1}
\cat(X)\geq \cwgt(z) +1
\end{eqnarray}
for any nonzero homology class $z\in H_q(X;R)$, $z\not=0$. The last inequality
can be also rewritten as
\begin{eqnarray}
0\leq \cwgt(z) \leq \cat(X)-1\leq \dim X\end{eqnarray} for any homology class.

Note that if $X$ is path-connected and $z$ is zero-dimensional, i.e. $z\in
H_0(X)$, then $\cwgt(z) = \cat(X)-1$.

\begin{lemma}
\label{coeffchange}
Let $f: R\to R'$ be a morphism of local coefficient systems over
$X$ and let $f_\ast: H_q(X;R)\to H_q(X;R')$ be the induced map on homology.
Then for any $z\in H_q(X;R)$ one has
\begin{eqnarray}
\cwgt(f_\ast(z)) \geq \cwgt(z).
\end{eqnarray}
\end{lemma}
\begin{proof}
The result follows directly by applying the definition.
\end{proof}

\begin{lemma} Assume that $X$ is a simplicial polyhedron. Then $\cwgt(z)\geq k$
iff $z$ can be realized in $X-A$ for any sub-polyhedron $A\subset X$ with
$\cat_XA\leq k$.
\end{lemma}
\begin{proof}
We only need to show the 'if'-direction. Let $A\subset X$ be
closed with $\cat_X A \leq k$. We need to show that $z$ can be
realized by a cycle in $X-A$. We have $A\subset U_1\cup\dots\cup
U_k$ with each $U_i$ open and null-homotopic in $X$. Passing to a
fine subdivision of $X$, we can find a sub-polyhedron $B\subset X$
with $A\subset B \subset U_1\cup \dots\cup U_k$. Then $\cat_XB\leq
k$ and $z$ can be realized by a cycle lying in $X-B\subset X-A$.
\end{proof}

\begin{example} Assume that $X$ is a closed 2-dimensional manifold, i.e. a
compact surface. Let us show that any nonzero homology class $z\in
H_1(X)$ has $\cwgt(z)\geq 1$. Indeed, it is easy to see that any
closed subset $A\subset X$ which is null-homotopic in $X$ lies in
the interior of a disk $D^2\subset X$; but $H_1(X- \Int D^2) \to
H_1(X)$ is an isomorphism.
\end{example}

\section{Homotopy invariance of $\cwgt(z)$}

\begin{lemma}
Let $f: X\to Y$ and $g: Y\to X$ be two continuous maps with $g\circ f\simeq
1_X$. Let $R'$ be a local coefficient system over $Y$ and $R=f^\ast R'$ be the
induced local system over $X$. Given a homology class $z\in H_q(X;R)$, define
$z'\in H_q(Y;R')$ by $z'=f_\ast(z)$. Then their category weights satisfy
\begin{eqnarray}
\cwgt(z')\geq \cwgt(z).
\end{eqnarray}
\end{lemma}
\begin{proof}
We start with the following well-known general remark. Let $B'\subset Y$ be a
subset which is null-homotopic in $Y$. Then the set $B=f^{-1}(B')\subset X$ is
null-homotopic in $X$. Indeed, since $1_X\simeq g\circ f$, the inclusion $B\to
X$ is homotopic to the composition $B\stackrel f\to B'\stackrel i\to Y\stackrel
g\to X$ where the inclusion $i: B'\to Y$ is null-homotopic by the assumption.

Denote $k=\cwgt(z)$. Assume that $A'\subset Y$ is a closed subset with
$\cat_YA'\leq k$. Consider $A=f^{-1}(A')\subset X$.  Since $\cat_YA'\leq k$
there exist open sets $U'_1, \dots, U'_k\subset Y$ covering $A'$ with each
$U'_i\to Y$ null-homotopic. Then the sets $U_i=f^{-1}(U'_i)\subset X$ are open,
cover $A$ and are null-homotopic in $X$ (by the above remark). This shows that
$\cat_XA\leq k$.

Since $\cwgt(z)\geq k$, the class $z$ can be realized by a singular cycle in
$X-A$. Then the cycle $c'=f_\ast(c)$ in $Y$ represents the class $z'$ and is
disjoint from $A'$ as $f$ maps $X-A$ into $Y-A'$.
\end{proof}

As a corollary of the previous result we obtain homotopy invariance of the
category weight:

\begin{theorem} If $f: X\to Y$ is a homotopy equivalence then for any homology
class $z\in H_q(X;R)$ one has
\begin{eqnarray}
\cwgt(z)=\cwgt(f_\ast(z)).
\end{eqnarray}
Here $f_\ast(z)\in H_q(X;R')$ where $R'=g^\ast R$ is the local coefficient
system over $Y$ induced by the homotopy inverse $g: Y\to X$ of $f$.
\end{theorem}

\section{Further properties of the category weight}

\begin{theorem}\label{cap1} Suppose that $X$ is a metric space.
Assume $u\in H^r(X;R)$, $z\in H_q(X;R')$ where $R$ and $R$ are local systems
over $X$. Then for the homology class $u\cap z\in H_{q-r}(X;R\otimes R')$ one
has
\begin{eqnarray}
\cwgt(u\cap z) \geq \cwgt(u) + \cwgt(z). \end{eqnarray} Here $\cwgt(z)$ is the
category weight of homology class $z$ as defined above in this paper and
$\cwgt(u)$ is the category weight of $u$ as defined by Fadell and Husseini
\cite{FH}.
\end{theorem}
\begin{proof} Denote $k=\cwgt(z)$, $\ell=\cwgt(u)$ and assume
that $A\subset X$ is a closed subset with $\cat_XA\leq k+\ell$. We want to show
that $u\cap z$ can be realized in the complement $X-A$. There exists an open
cover $A\subset U_1\cup U_2\cup \dots\cup U_{k+\ell}\subset X$ with each
$U_i\to X$ null-homotopic. Find open subsets $V_i\subset U_i$ such that $\bar
V_i\subset U_i$ and $A\subset V_1\cup \dots\cup V_{k+\ell}$.

Denote $B= V_1\cup\dots\cup V_\ell$ and let $C=A-B$. Clearly $C$ is closed and
satisfies $\cat_XC\leq k$ Hence $z$ can be realized by a cycle avoiding $C$. In
other words, $z=i_\ast(w)$ where $w\in H_q(X-C;R')$.

Since $\cwgt(u)\geq \ell$ we have $u|_{\overline B}=0$ and thus $u=j_1^\ast(v)$
for some $v\in H^r(X,B;R)$. By statements 16 in \cite{Sp}, chapter 5, \S 6, one
has
$$j_\ast(u\cap z) = j_\ast(j_1^\ast v\cap z) = v \cap \bar j_\ast (i_\ast w ) =0$$
where $j: X\to (X, X-A)$, $\bar j: X \to (X, X-C)$ and $j_1: X\to (X,B)$ are
inclusions. By exactness, $j_\ast(u\cap z)=0$ implies that $u\cap z$ lies in
the image of $H_{q-r}(X-A;R\otimes R') \to H_{q-r}(X;R\otimes R')$.
\end{proof}

As a corollary we obtain:

\begin{corollary} Suppose that $X$ is a metric space and for some classes
$z\in H_q(X;R)$ and $u\in H^q(X;R')$ the evaluation
$$\langle u, z\rangle \, \not=0\, \in\,  R'\otimes R$$ is nonzero.
Then
\begin{eqnarray}\label{improved}
\cat(X) \geq \cwgt(z) +\cwgt(u) +1.
\end{eqnarray}
Here $\cwgt(z)$ is the category weight of homology class $z$ as defined above
in this paper and $\cwgt(u)$ is the category weight of $u$ as defined by Fadell
and Husseini \cite{FH}.
\end{corollary}

\begin{proof} This follows from inequality (\ref{one1}) combined with Theorem
\ref{cap1}.
\end{proof}
Inequality (\ref{improved}) allows us to improve the classical cohomological lower
bound for the category $\cat(X)$ by taking into account the quality of the
homology class $z$.

\section{Manifolds and Poincar\'e complexes}

In this section we prove that in the case of closed manifolds our notion of
category weight coincides with the cohomological notion of Fadell and Husseini
\cite{FH}. However for Poincar\'e complexes these notions are distinct as we
show by an example.

\begin{theorem}\label{manif}
Suppose that $X$ is a closed $n$-dimensional manifold, $z\in H_q(X;R)$ where
$R$ is a local coefficient system. Let $u\in H^{n-q}(X; R\otimes \tilde \Z)$ be
the Poincar\'e dual cohomology class, i.e. $z=u\cap [X]$, see below. Then
\begin{eqnarray} \cwgt(z)
=\cwgt(u).\end{eqnarray} Here $\tilde \Z$ denotes the orientation local system
on $X$, i.e. for a point $x\in X$ the stalk of $\tilde \Z$ at $x$ is $\tilde
\Z_x=H_n(X, X-x;\Z)$, see \cite{Sp}.
\end{theorem}

\begin{proof}
By Poincar\'e duality theorem any homology class $z\in H_q(X;R)$ can be
uniquely written as $z=u\cap [X]$ where $u\in H^{n-q}(X;R\otimes \tilde \Z)$
and $[X]\in H_n(X;\tilde \Z)$ is the fundamental class.  Applying inequality of
Theorem \ref{cap1} we find
\begin{eqnarray}\label{inequality}
\cwgt(z) \geq \cwgt(u)+\cwgt([X]) =
\cwgt(u).
\end{eqnarray}
To obtain the inverse inequality one observes that if
$A\subset X$ is a closed subset with $\cat_X A\leq \cwgt(z)$ then $z$ can be
realized by a singular cycle in the complement $X-A$ and the usual intersection
theory for chains in manifolds shows that the cocycle Poincar\'e dual to $z$
vanishes on $A$; hence $\cwgt(u)\geq \cwgt(z)$.
\end{proof}

\begin{example}
Let $X=\RP^n$ be the real projective space. For the unique nonzero
cohomology class $z\in H_q(X;\Z_2)$ one has $\cwgt(z) = n-q.$
Indeed, the dual homology class is $\alpha^{n-q}\in
H^{n-q}(X;\Z_2)$ where $\alpha\in H^1(X;\Z_2)$ is the generator.
Clearly, $\cwgt(\alpha^{n-q}) = n-q$.
\end{example}

Theorem \ref{manif} implies:

\begin{corollary}
If $X$ is a closed $n$-dimensional manifold then for any homology class $z\in
H_q(X;R)$ with $q<n$ one has \begin{eqnarray} \cwgt(z)\geq 1.\end{eqnarray}
\end{corollary}

Indeed, if $q<n$ then the dual cohomology class $u$ has positive degree and
hence $\cwgt(u)\geq 1$.

Consider now the case when $X$ is $n$-dimensional Poincar\'e
complex. The first part of proof of Theorem \ref{manif} is still
applicable giving inequality (\ref{inequality}) between category
weights of the homology and cohomology classes. However the second
part of the proof fails. The following example shows that Theorem
\ref{manif} is false for Poincar\'e complexes. It is a
modification of an argument due to D. Puppe showing that the
notion of category weight of cohomology classes is not homotopy
invariant.

\begin{example} Consider the lens space $L=S^{2n+1}/(\Z/p)$ where $p$ is
an odd prime and $\Z/p$ acts freely on $S^{2n+1}$. Denote by $r: S^{2n+1}\to L$
the quotient map. Let $X$ be the mapping cylinder of $r$, i.e.
$$X= L \sqcup S^{2n+1}\times [0,1]/\sim$$
where each point $(x, 0)\in S^{2n+1}\times [0,1]$ is identified with $r(x)\in
L$. Clearly $X$ is homotopy equivalent to $L$ and so it is a Poincar\'e
complex. By a theorem of Krasnoselski \cite{Kr}, the category of $X$ equals
$2n+2$. Hence for $z=1\in H_0(X;\Z_2)$ one has
$$\cwgt(z)=\cat(X)-1=2n+1,$$ see
above. The dual cohomology class $u$ is the generator $u\in
H^{2n+1}(X;\Z_2)$. Let us show that
$$\cwgt(u)=1.$$ Indeed, consider the sphere
$S=S^{2n+1}\times 1\subset X$. The restriction $u|_S\in H^{2n+1}(S;\Z_2)$
coincides with the induced class $r^\ast(v)$ where $v\in H^{2n+1}(L;\Z_2)$ is
the generator. Hence the cohomology class $u|_S$ is nonzero. However, the
sphere $S$ has category 2 and moreover $\cat_XS=2$ (as the inclusion $S\to X$
is not null-homotopic).
\end{example}

The following simple construction gives non-manifolds for which the category
weight can be explicitly calculated.

\begin{lemma}
Let $X=X_1\vee X_2$ be wedge of two polyhera $X_1$ and $X_2$ and let $z\in
H_q(X;R)$ be the sum $z=z_1+z_2$ where $z_i\in H_q(X_i; R_i)$ and $R_i
=R|_{X_i}.$ Then
\begin{eqnarray}
\cwgt(z) = \min\{\cwgt(z_1), \cwgt(z_2)\}.
\end{eqnarray}
Here $\cwgt(z_i)$ is the category weight of $z_i$ viewed as a homology class of
$X_i$.
\end{lemma}
\begin{proof}
Inequality $\cwgt(z) \leq \min\{\cwgt(z_1), \cwgt(z_2)\}$ is obvious. Let
$A\subset X$ be a closed subset with $\cat_XA\leq k$ where $k=
\min\{\cwgt(z_1), \cwgt(z_2)\}$. Then $A=A_1\vee A_2$ where $A_i\subset X_i$
and $\cat_{X_i}A_i\leq k$, where $i=1, 2$. One can realize $z_i$ by a cycle
avoiding $A_i$. The sum of these two cycles is a cycle representing $z$ which
avoids $A$. Thus we obtain the opposite inequality $\cwgt(z)\geq k$.
\end{proof}

\section{Strict category weight}

The notion of strict category weight was introduced in \cite{Ru};
it is a homotopy invariant variation of the category weight of
Fadell and Husseini \cite{FH}. We use this notion in this paper
and therefore recall the relevant definitions. We warn the reader
that our terminology differs from \cite{Ru} by 1 and is consistent
with \cite{FH}.

\begin{definition}
Given a continuous map $\phi: A\to X$, we say that $\cat(\phi) \leq k$ if $A$
can be covered by $k$ open sets $A_1, \dots, A_k$ such that each restriction
$\phi|_{A_i}$ is null-homotopic. The strict category weight of a cohomology
class $u\in H^q(X;R)$ (where $R$ is a local coefficient system on $X$) is
defined as the maximal integer $k$ such that $\phi^\ast(u)=0$ for any
continuous map $\phi: A\to X$ with $\cat(\phi)\leq k$.
\end{definition}

The strict category weight is denoted by $\swgt (u)$. Clearly, one has
$$\swgt(u)\leq \cwgt(u)$$ and $\swgt(u)\geq 1$ for any cohomology class $u\in
H^q(X;R)$ of positive degree $q>0$.

\begin{definition}
Let $X$ be a closed smooth connected $n$-dimensional manifold. We define the
strict category weight of a homology class $z\in H_q(X;R)$ (denoted $\swgt(z)$)
as the strict category weight of the dual cohomology class $u\in
H^{n-q}(X;R\otimes \tilde Z)$.
\end{definition}

Similar definition can be used in the case of Poincar\'e complexes, but we do
not use it in such generality.

\begin{proposition}\label{strict}
Let $z_i\in H_{q_i}(X_i, R_i)$ where $X_i$ is a closed smooth orientable
manifold of dimension $n_i$, $i=1, 2$. Consider the cross-product
$$z_1\times z_2\in H_{q}(X_1\times X_2; R)$$ where $q=q_1+q_2$ and $R$ is the
external tensor product $R=R_1\boxtimes R_2$. Then
\begin{eqnarray}
\swgt(z_1\times z_2) \geq \swgt(z_1) +\swgt(z_2).
\end{eqnarray}
\end{proposition}
\begin{proof}
Let $u_i\in H^{n_i-q_i}(X_i;R_i)$ denote the dual of $z_i$, where $i=1,2$. Then
the dual of $z_1\times z_2$ is $u_1\times u_1\in H^{n-q}(X_1\times X_2;R)$
where $n=n_1+n_2$. Consider also the classes $u_1\times 1\in
H^{n_1-q_1}(X_1\times X_2;R_1\boxtimes \Z)$ and $1\times u_2\in
H^{n_2-q_2}(X_1\times X_2;\Z\boxtimes R_2)$.

Denote $k_i=\swgt(z_i)=\swgt(u_i)$. Let $\phi: A\to X_1\times X_2$ be a
continuous map with $\cat \phi \leq k_1+k_2$. Then $A$ is union of open subsets
$A=A_1\cup A_2$ such that $\cat \phi|_{A_i}\leq k_i$. We obtain that
$\phi^\ast(u_1\times 1)|_{A_1}=0$ and $\phi^\ast(1\times u_2)|_{A_2}=0$. This
implies that the class $\phi^\ast(u_1\times u_2) = \phi^\ast(u_1\times 1)\cup
\phi^\ast(1\times u_2)$ vanishes. Hence $\swgt(z_1\times z_2) \geq k_1+k_2$.
\end{proof}

\begin{corollary}\label{prodweight}
Let $X_i$ be closed orientable manifolds and $z_i\in H_{q_i}(X_i;R_i)$ where
$q_i<\dim X_i$ for $i=1, \dots, k$. Consider $z=z_1\times \dots\times z_k\in
H_q(X;R)$ where $X=X_1\times \dots\times X_k$, $q=q_1+\dots+q_k$ and
$R=R_1\boxtimes \dots\boxtimes R_k$. Then
\begin{eqnarray}
\cwgt(z)  \geq k.
\end{eqnarray}
\end{corollary}

This corollary is a source of examples of homology classes having
high category weight.

\specialsection*{\bf Part II: Moving integral homology classes to infinity}

In Part II we study conditions for an integral homology class
$z\in H_i(\tilde X;\Z)$ of a free abelian covering $\tilde X\to X$
to be movable to infinity with respect to a cohomology class
$\xi\in H^1(X;\R)$. The case of homology classes with coefficients
in a field was studied in \cite{FS} using a different algebraic
technique which is not applicable over $\Z$.

\section{Abel - Jacobi maps and neighborhoods of infinity}

\label{sec:aj}

For the convenience of the reader we recall in this section the
language introduced in \cite{FS}. Let $X$ be a connected finite
cell complex and $p:\tilde X\to X$ a regular covering having a
free abelian group of covering transformations $H\simeq \Z^r$.
Denote $ H_\R=H\otimes \R; $ it is a vector space of dimension $r$
containing $H$ as a lattice.
\begin{proposition}
There exists a canonical Abel - Jacobi map
\begin{eqnarray} A: \tilde X\to H_\R
\end{eqnarray}
having the following properties:
\begin{enumerate}
\item[(a)] $A$ is $H$-equivariant; here $H$ acts on $\tilde X$ by covering
transformations and it acts on $H_\R$ by translations.

 \item[(b)] $A$ is proper (i.e. the preimage of
a compact subset of $H_\R$ is compact).

\item[(c)] $A$ is determined uniquely up to replacing it by a map $A': \tilde
X\to H_\R$ of the form $ A' = A+ F\circ p $ where $F: X\to H_\R$ is a
continuous map.
\end{enumerate}\end{proposition}

This fact is well-known; we refer to \cite{FS} for a detailed proof.

Let $\xi\in H^1(X;\R)$ be a cohomology class with the property
$$p^\ast(\xi)=0\in H^1(\tilde X;\R).$$ Such a class $\xi$ can be viewed either as a
homomorphism $\xi: H\to \R$ or as a linear functional $\xi_\R: H_\R\to \R$.

\begin{definition}
A subset $N\subset \tilde X$ is called a neighborhood of infinity in $\tilde X$
with respect to the cohomology class $\xi$ if $N$ contains the set
\begin{eqnarray}
\{x\in \tilde X; \, \, \xi_\R(A(x))>c\}\, \subset \, N,
\end{eqnarray}
for some real $c\in \R$. Here $A: \tilde X\to H_\R$ is an Abel - Jacobi map for
the covering $p:\tilde X\to X$.
\end{definition}

See \cite{FS} for more details.

\section{Homology classes movable to infinity}

Let $G$ be an abelian group (the coefficient system). We mainly have in mind
the cases $G=\Z$ or $G=\kk$ is a field.

\begin{definition} (See \cite{farbe2}, \S 5)
A homology class $z\in H_i(\tilde X;G)$ is said to be movable to infinity of
$\tilde X$ with respect to a nonzero cohomology class $\xi\in H^1(X;\R)$,
$p^\ast(\xi)=0$, if in any neighborhood $N$ of infinity with respect to $\xi$
there exists a (singular) cycle with coefficients in $G$ representing $z$.
\end{definition}

Equivalently, a homology class $z\in H_i(\tilde X;G)$ is movable to infinity
with respect to $\xi\in H^1(X;\R)$ if $z$ lies in the intersection
\begin{eqnarray}
\bigcap_N  \im\left[ H_i(N;G)\to H_i(\tilde X;G)\right]
\end{eqnarray}
where $N$ runs over all neighborhoods of infinity in $\tilde X$ with respect to
$\xi$. This can also be expressed by saying that $z$ lies in the kernel of the
natural homomorphism
\begin{eqnarray}\label{kernel}
H_i(\tilde X;G) \to \lim_\leftarrow H_i(\tilde X,N;G)
\end{eqnarray}
where in the inverse limit $N$ runs over all neighborhoods of infinity in
$\tilde X$ with respect to $\xi$.


The following theorem proven in \cite{FS} gives an explicit description of all
movable homology classes in the case when $G=\kk$ is a field. It generalizes
the result of \cite{farbe2}, \S 5 treating the simplest case of infinite cyclic
covers $q: \tilde X\to X$.

\begin{theorem}\label{field} Let $X$ be a finite cell complex and $q:\tilde X\to X$
be a regular covering having a free abelian group of covering transformations
$H\simeq\Z^r$. Let $\xi\in H^1(X;\R)$ be a nonzero cohomology class of rank $r$
satisfying $q^\ast(\xi)=0$. The following properties (A), (B), (C) of a nonzero
homology class $z\in H_i(\tilde X;\kk)$ (where $\kk$ is a field) are
equivalent:

\begin{enumerate}
\item[(A)] $z$ is movable to infinity with respect to $\xi$.

\item[(B)] Any singular cycle $c$ in $\tilde X$ realizing the class $z$ bounds
an infinite singular chain $c'$ in $\tilde X$ containing only finitely many
simplices lying outside every neighborhood of infinity $N\subset \tilde X$ with
respect to $\xi$.

\item[(C)] There exists a nonzero element $x\in \kk[H]$ such that $x\cdot z=0$.
\end{enumerate}
\end{theorem}

Later in this paper (see \S \ref{sec7}) we will describe the set of homology
classes with integral coefficients which are movable to infinity.

\section{Integral homology classes movable to infinity}\label{sec7}

To get an analogue of Theorem \ref{field} in the case of integral coefficients, we need another definition.

\begin{definition}
Let $H$ be a group and $\xi:H\to \R$ a homomorphism. A nonzero element $\Delta\in \Z[H]$ is said to have $\xi$-lowest coefficient 1, if $\Delta=(1-y)h$ with $h\in H$ and $y=\sum a_jg_j$, where the $g_j\in H$ satisfy $\xi(g_j)>0$ and $a_j\in \Z$.
\end{definition}

\begin{theorem}\label{movability} Let $X$ be a finite cell complex and $p:\tilde X\to X$
be a regular covering having a free abelian group of covering transformations
$H\simeq\Z^r$. Let $\xi\in H^1(X;\R)$ be a nonzero cohomology class of rank $r$
satisfying $p^\ast(\xi)=0$. The following properties (A), (B), (C) of a nonzero
integral homology class $z\in H_i(\tilde X;\Z)$ are equivalent:

\begin{enumerate}
\item[(A)] $z$ is movable to infinity with respect to $\xi$.

\item[(B)] Any singular cycle $c$ in $\tilde X$ realizing the class $z$ bounds
an infinite singular chain $c'$ in $\tilde X$ with integral coefficients
containing only finitely many simplices lying outside every neighborhood of
infinity $N\subset \tilde X$ with respect to $\xi$.

\item[(C)] There exists a nonzero element $\Delta\in \Z[H]$ with $\xi$-lowest
coefficient 1 such that $\Delta\cdot z=0$ .
\end{enumerate}
\end{theorem}

 This result improves Theorem 5.3 of \cite{farbe4} which treats the case of
 rank one cohomology classes, $r=1$. Movability to infinity of homology classes
 with coefficients in a field was studied in \cite{farbe2} ($r=1$ case) and in \cite{FS}
 ($r\geq 1$).

Note that the implications (C) $\Rightarrow$ (B) $\Rightarrow$ (A) of Theorem
\ref{movability} are straightforward (see below); the only nontrivial statement
is the implication (A) $\Rightarrow $ (C). Let us explain why (C) $\Rightarrow$
(B). Suppose that $\Delta \cdot z =0\in H_i(\tilde X;\Z)$ where $\Delta\in
\Z[H]$ has $\xi$-lowest coefficient 1. Without loss of generality we may assume
that $\Delta=1-y$ where $y\in \Z[H]$ is $\xi$-positive, i.e. $y$ is a finite
sum of the form $\sum a_j g_j$ where $g_j\in H$, $\xi(g_j)>0$, and $a_j \in
\Z$. Let $c$ be a chain representing the class $z$. Then the cycle $\Delta\cdot
c$ bounds, i.e. $(1-y)\cdot c=\partial c_1$ where $c_1$ is a finite chain in
$\tilde X$. Set $c'= c_1 + yc_1+y^2c_1+\dots.$ Then $\partial c'=c$ and $c'$
has finitely many simplices lying outside every neighborhood of infinity
$N\subset \tilde X$ with respect to $\xi$.

The main part of the proof consists in establishing the vanishing of the
$\lim^1$ term in the following exact sequence:

\begin{eqnarray}\label{exact1}\quad \quad 0\to
{\lim_{\leftarrow}}^1 H_{q+1}(\tilde X,N;\Z)\to H_q(X;\widehat{\Z[H]}_\xi)
\to {\lim_{\leftarrow}}\,  H_q(\tilde X,N;\Z)\to 0.
\end{eqnarray}
This exact sequence was described in \S 6 of \cite{FS}. Formally,
the proof of the exactness (\ref{exact1}) given in \cite{FS}
assumes that the ring of coefficients is a field but it works
equally well in the case $\Z$ with no modifications. In the exact
sequence (\ref{exact1}) $\lim$ and $\lim^1$ are taken relative to
the system of neighborhoods of infinity $N\subset \tilde X$ with
respect to $\xi$. The symbol $\widehat{\Z[H]}_\xi$ in
(\ref{exact1}) denotes the Novikov completion of the group ring
$\Z[H]$, see \cite{noviko1}, \cite{noviko2}. Recall that elements
of the group ring $\Z[H]$ are finite sums of the form $ \sum
a_ig_i$ where $a_i\in \Z$ and $g_i\in H$; the ring
$\widehat{\Z[H]}_\xi$ includes also all countable sums $\sum
a_ig_i$ having the property $ \lim_{i\to +\infty}\xi(g_i)\, =\,
+\infty. $

\begin{proposition}\label{prop3}
Under the conditions of Theorem \ref{movability} one has
\begin{eqnarray}\label{vanish}
{\lim_{\leftarrow}}^1 H_{q}(\tilde X,N;\Z)=0,\end{eqnarray} where
$N$ runs over all neighborhoods of infinity in $\tilde X$ with
respect to $\xi$ partially ordered by the inverse inclusion.
\end{proposition}

Proposition \ref{prop3} gives the implication (A) $\Rightarrow$ (B) of Theorem
\ref{movability}. Indeed, using the definition (\ref{kernel}) combined with
(\ref{vanish}) we see that a homology class $z\in H_q(\tilde X;\Z)$ is movable
to infinity with respect to $\xi$ if and only if a cycle $c\in C_q(\tilde X)$
representing $z$ bounds a chain $c'\in C_q(\tilde X)\otimes
\widehat{\Z[H]}_\xi$, i.e. $dc'=c$. Here $C_\ast(\tilde X)$ denote the cellular
chain complex of $\tilde X$ with integral coefficients. One can view $c'$ as an
infinite chain in $\tilde X$ having finitely many terms outside any given
neighborhood of infinity in $\tilde X$ with respect to $\xi$.

To see that (B) $\Rightarrow$ (C), let $S_\xi\subset \Z[H]$ be the subset consisting of elements with $\xi$-lowest coefficient 1 and $\Lambda_\xi=S^{-1}_\xi \Z[H]$ the localization. By \cite[Lemma 1.13]{farbe3} the inclusion $\Lambda_\xi\to \widehat{\Z[H]}_\xi$ is faithfully flat so that the change of coefficients $H_\ast(X;\Lambda_\xi)\to H_\ast(X;\widehat{\Z[H]}_\xi)$ is injective as well. The result follows.

\section{Proof of Theorem \ref{movability}.}

First we discuss some commutative algebra.

Recall our notations. $H=\Z^r$ is a free abelian group and $\xi: H\to \R$ is an
injective group homomorphism. We denote by $A$ the Novikov ring
$\widehat{\Z[H]_\xi}$ and by $A_0$ its subring $\widehat{\Z[H_0]_\xi}$ where
$H_0=\{g\in H; \xi(g)\geq 0\}$. Elements of $A_0$ are countable formal sums of
the form $\sum_j a_j g_j$ where $a_j\in \Z$ and $\xi(g_j)$ tends to $+\infty$.

It is well-known that $A$ is a principal ideal domain but $A_0$ is not. Our
goal is to obtain some partial results about properties of modules over the
ring $A_0$ resembling those of modules over principle ideal domains.

\begin{definition}
Let $M$ be a $A_0$-module. A sequence of elements $m_1, \dots, m_k\in M$ is a
quasi-basis for $M$ if (1) for any $m\in M$ there exists $g\in H_0$ such that
$gm$ can be represented in the form $gm = \sum a_j m_j$ where $a_j\in A_0$, and
(2) there are no nontrivial relations $\sum a_j m_j=0$.
\end{definition}

\begin{lemma}\label{quasi}
Let $f: A_0^n \to A_0^m$ be a homomorphism of finitely generated free
$A_0$-modules. Then there exist quasi-basis $d_1, \dots, d_n\in A_0^n$ and
$e_1, \dots,  e_m\in A_0^m$ and an integer $\mu \leq \min\{n, m\}$ such that
for any $j\leq \mu$
 one has
\begin{eqnarray}\label{relations}
f(d_j) = a_j e_j, \quad \mbox{where}\quad a_j\in A_0, \quad a_j\not=0,
\end{eqnarray}
and $ f(d_j)=0 $ for $j>\mu$.
\end{lemma}
\begin{proof}
Localizations of $A_0^n$ and $A_0^m$ with respect to the multiplicative set $H$
lead to free modules over principles ideal domain $A$. Hence, applying the
standard theory, we find free basis $d_1', \dots, d_n'\in A^n$, $e_1', \dots,
e'_m$ and an integer $\mu\leq \max\{n, m\}$ such that $f(d'_j)=a_j'e'_j$ for
$j\leq \mu$ (where $a'_j\in A$, $a'_j\not=0$) and $f(d'_j)=0$ for $j>\mu$.
Choose $g\in H_0$ such that $gd'_j\in A_0^n$ and $ge'_j\in A_0^m$ for all $j$.
Choose $g'\in H_0$ such that $g'a'_j\in A_0$ for all $j$. Now set $d_j =
gg'd'_j$, $e_j=ge'_j$, and $a_j =g'a'_j$. We obtain quasi-basis $d_j$ and $e_j$
and clearly the relations (\ref{relations}) hold.
\end{proof}
\begin{lemma}\label{subcomplex}
Let $C_\ast$ be a free finitely generated\footnote{By this we mean that each
$A_0$-module $C_q$ is finitely generated and only finitely many modules $C_q$
are nonzero.} chain complex over $A_0$. Then there exist quasi-basis $e^q_1,
e^q_2, \dots, e^q_{r_q}\in C_q$ (where $q\in \Z$ and $r_q$ denotes the rank of
$C_q$) and integers $\mu_q\leq \min\{r_q, r_{q-1}\}$ such that the differential
$d: C_q\to C_{q-1}$ is given by
\begin{eqnarray}
d(e^q_j) = \left\{
\begin{array}{lll}
a^q_j e^{q-1}_j & \mbox{for}& j\leq \mu_q,\\
0 & \mbox{for}& j > \mu_q,
\end{array}
\right.
\end{eqnarray}
and the elements $a^q_j\in A_0$ are nonzero.
\end{lemma}
\begin{proof}
The proof essentially repeats the arguments of Lemma \ref{relations}. On the
first step we construct basis $f^q_j$ of the localized chain complex
$C'_q=A\otimes_{A_0}C_q$ over the principal ideal domain $A$ such that all
differentials $d: C'_q\to C'_{q-1}$ have the diagonal form $d(f^q_j) =
\alpha^q_j f^{q-1}_j$ with $\alpha^q_j\in A$. On the second step one multiplies
the basis $f^q_j$ by a suitable group element $g^q\in H_0$ so that (1) the
elements $e^q_j=g^q f^q_j$ lies in the original complex $C_q$ and (2) the
elements $a^q_j = g^q (g^{q-1})^{-1} \alpha^q_j$ lie in $A_0$.

\end{proof}

\begin{lemma}\label{submodule}
Let $C_\ast$ be a free finitely generated chain complex over $A_0$. Then there
exists a finitely generated free chain subcomplex $D_\ast\subset C_\ast$ such
that $gC_\ast\subset D_\ast$ for some group element $g\in H_0$ and $
H_j(D_\ast) $ is isomorphic to a finite direct sum of cyclic\footnote{i.e.
modules of the form $A_0/(aA_0)$ where $a\in A_0$.} $A_0$-modules.
\end{lemma}
\begin{proof} Apply Lemma \ref{subcomplex} and take for $D_q\subset C_q$ the
$A_0$-submodule generated by $e^q_1, \dots, e^q_{r_q}$.
\end{proof}

Next we apply the above results to obtain the following corollary.

\begin{corollary}\label{cor3}
Let $C_\ast$ be a free finitely generated chain complex over $A_0$. Let
$\bar{C_\ast}= A\otimes_{A_0}C_\ast$  be the localized chain complex and $i:
C_\ast \to \bar{C_\ast}$ the inclusion. Then for any $q$ there exists a group
element $g=g^q\in H_0$ such that the kernel of the induced map $i_\ast: H_q(C)
\to H_q(\bar{C})$ coincides with the kernel of multiplication by $g$ on
$H_q(C)$.
\end{corollary}
\begin{proof}
As a preparation, consider a nonzero cyclic $A_0$-module $M_0=A_0/(aA_0)$ and
the associated $A$-module $M=A/(aA)$. Here $a\in A_0$ is a non-invertible
element. Write $a$ in the form $a=g(\alpha + h\beta)$ where $g\in H_0$,
$\alpha\in \Z$, $\alpha\not=0$, $\beta\in A_0$ and $h\in H$ is such that
$\xi(h)>0$. Note that $M$ is trivial if and only if $a$ is invertible in $A$,
i.e. when $\alpha=\pm 1$. Similarly, $M_0$ is trivial iff $a$ is invertible in
$A_0$ i.e. when $g=0$ and $\alpha=\pm 1$. We will say that $M_0$ is a cyclic
module of the first (second) kind if $\alpha=\pm 1$ (or $|\alpha|>1$,
correspondingly). We obtain that for a cyclic module $M_0$ of the first kind
there exists $g\in H_0$ such that $gM_0=0$ and the corresponding module $M$ is
trivial. For a cyclic module $M_0$ of the second kind, there is a $g\in H$ such that $gM_0\to M$ is injective. Indeed, with $a=g(\alpha+h\beta)$ as above and $M_0=A_0/(aA_0)$, we get $gM_0=A_0/(\alpha+h\beta A_0)$ which includes into $A/(\alpha+h\beta A)$.

Apply Lemma \ref{submodule} to obtain a subcomplex $D_\ast\subset C_\ast$ such
that $g'C_\ast\subset D_\ast$ for some $g'\in H_0$ and $H_\ast(D)$ is a finite
direct sum of cyclic $A_0$-modules. One finds that $A\otimes_{A_0}
D_\ast=A\otimes_{A_0} C_\ast =\bar C_\ast$ and $H_q(\bar C) = A\otimes_{A_0}
H_q(C) = A\otimes_{A_0} H_q(D)$ since $A$ is flat over $A_0$.
Also, the kernel of the map restricted to the summands of cyclic modules of the second kind can be annihilated by multiplication with a suitable element of $H_0$.
Thus there exists an $h\in H_0$ such that the kernel of the map $H_q(D)\to H_q(\bar{C})$ coincides with the kernel of the map $h:H_q(D)\to H_q(D)$.

Now set $g=hg'\in H_0$. Let us show that the kernel of $i_\ast: H_q(C) \to
H_q(\bar C)$ coincides with the kernel of multiplication $g_\ast: H_q(C)\to
H_q(C)$ by $g$. Consider the following diagram
\begin{eqnarray*}
\begin{array}{ccccclr}
H_q(C) & \stackrel {g'_\ast}\to & H_q(D)& \stackrel {h_\ast}\to & H_q(D)& \to
& H_q(C)\\
i_\ast \downarrow && j_\ast\downarrow && j_\ast \downarrow & & \downarrow i_\ast \\
H_q(\bar C) & \stackrel {g'_\ast}\to & H_q\bar C) & \stackrel {h_\ast}\to &
H_q(\bar C) &\stackrel {\rm {id}}\to & H_q(\bar C)
\end{array}
\end{eqnarray*}
The composition of the upper horizontal row is the multiplication by $g$, i.e.
the map $g_\ast: H_q(C)\to H_q(C)$. Every map appearing in the lower horizontal
row is an isomorphism. From the previous paragraph we know that $\ker (j_\ast)
= \ker (h_\ast)$. Therefore, examining the diagram, we find that $\ker(i_\ast)
= \ker (g_\ast)$ as claimed.
\end{proof}

The vanishing of the $\lim^1$ term of the exact sequence (\ref{exact1}) (i.e.
Proposition \ref{prop3}, see above) would follow once one has the Mittag-Lefler
condition (see \cite{W}, Prop. 3.5.7) which in our case states:

\begin{proposition}\label{mittag2}
For any neighborhood of infinity $N\subset \tilde X$ with respect to $\xi$
there exists a neighborhood of infinity $N'\subset N$ such that for any
neighborhood of infinity $N''\subset N'$ one has
 \begin{eqnarray}\label{mittag}
\quad \im \left[H_q(\tilde X,{N}'') \to H_q(\tilde X,N)\right] =\im
\left[H_q(\tilde X,N') \to H_q(\tilde X,N)\right].
 \end{eqnarray}
 \end{proposition}
The homology groups appearing in
 (\ref{mittag})
are with coefficients in $\Z$ and all neighborhoods of infinity are with
respect to a fixed cohomology class $\xi$.

The equality (\ref{mittag}) can be expressed by saying that any cycle in $\tilde X$
relative to $N$ which can be refined to a cycle relative to $N'$ can be refined
to a cycle relative to an arbitrary neighborhood of infinity $N''\subset
N'\subset \tilde X$.

\begin{proof}[Proof of Proposition \ref{mittag2}.]
Let $C_\ast(\tilde X)$ denote the cellular chain complex of $\tilde X$. It is a
complex of finitely generated free $\Z[H]$-modules. Let $N\subset \tilde X$ be
a cellular neighborhood of infinity with respect to $\xi$ as described in Lemma
3 of \cite{FS}. The cellular chain complex $C_\ast(N)$ is free and finitely
generated over $\Z[H_0]$ where $H_0=\{g\in H; \xi(g)\geq 0\}$. Consider the
completed chain complexes $C'_\ast(N)= A_0\otimes_{\Z[H_0]}C_\ast(N)$ and
$C'_\ast(\tilde X) = A\otimes_{\Z[H]}C_\ast(\tilde X)$. Recall that
$A=\widehat{\Z[H]_\xi}$ is the Novikov ring and $A_0=\widehat{\Z[H_0]_\xi}$.
The canonical inclusions $C_\ast(N)\to C'_\ast(N)$ and $C_\ast(\tilde X)\to
C'_\ast(\tilde X)$ determine a chain homomorphism
\begin{eqnarray}\label{isom}
F:C_\ast(\tilde X)/C_\ast(N)\stackrel \simeq\to C'_\ast(\tilde X)/C'_\ast(N)
\end{eqnarray}
which is an isomorphism. Injectivity of $F$ is equivalent to $C_\ast(\tilde
X)\cap C'_\ast(N)=C_\ast(N)$ (which is obvious) and surjectivity of $F$ is
equivalent to $C_\ast(\tilde X)+C'_\ast(N)=C'_\ast(\tilde X)$. The latter
follows from the following equality $\Z[H]+A_0=A$ for subrings of $A$.

The short exact sequence of chain complexes over $A_0$
$$0\to C'_\ast(N)\to C'_\ast(\tilde X) \to C_\ast(\tilde X)/C_\ast(N)\to 0$$
gives the exact sequence
$$\dots \to H'_q(N) \stackrel {i_\ast}\to H'_q(\tilde X) \to H_q(\tilde X, N)\stackrel \partial
\to H'_{q-1}(N) \to \dots$$ where $H'_\ast(N)$ denotes homology of the complex
$C'_\ast(N)$ and similarly for $H'_\ast(\tilde X)$; the symbol $H_q(\tilde X,
N)$ denotes $H_q(\tilde X, N;\Z)$.

Applying Corollary \ref{cor3} to the subcomplex $C'_\ast(N)\subset
C'_\ast(\tilde X)$ we find a group element $g\in H_0$ such that $\ker[i_\ast:
H'_{q-1}(N)\to H'_{q-1}(\tilde X)]$ coincides with $\ker[g_\ast: H'_{q-1}(N)
\to H'_{q-1}(N)]= \ker[j_\ast: H'_{q-1}(N)\to H'_{q-1}(g^{-1}N)]$. Here $j:
N\to g^{-1}N$ is the inclusion. Denoting $N'=gN\subset N$ we obtain:
\begin{eqnarray}\label{obtain}
\qquad \ker[i_\ast: H'_{q-1}(N')\to H'_{q-1}(\tilde X)] = \ker[i_\ast:
H'_{q-1}(N')\to H'_{q-1}(N)].
\end{eqnarray}

Now, consider the following commutative diagram
 \begin{eqnarray*}
 \begin{array}{ccccc}
 H'_q(N) & \to & H'_q(\tilde X) & \stackrel \beta \to & H_q(\tilde X, N)\\
 \downarrow && \downarrow && \downarrow =\\
 H_q(N, N') & \stackrel \tau \to & H_q(\tilde X, N') & \stackrel \alpha \to &
 H_q(\tilde X, N)\\
 \downarrow && \downarrow \partial &&\\
 H'_{q-1}(N') & \underset{=}{\stackrel \sigma \to}& H'_{q-1}(N') &&\\
 \downarrow \gamma &&\downarrow&&\\
 H'_{q-1}(N) &\to& H'_{q-1}(\tilde X)&&
 \end{array}
 \end{eqnarray*}
Clearly $\im \beta \subset \im \alpha$. The inverse implication $\im \alpha
\subset \im \beta$ would follow once we know that for any $x\in H_q(\tilde X,
N')$ there exists $y\in H_q(N, N')$ such that $\partial \tau (y) =\partial
(x)\in H'_{q-1}(N')$. Now, equality (\ref{obtain}) says that $\gamma \circ
\sigma^{-1}\circ
\partial=0$ is trivial which (using exactness properties of the diagram above)
means that for any $x\in H_q(\tilde X, N')$ an element $y\in H_q(N, N')$ with
the above mentioned property exists. This shows that $\im \alpha=\im \beta$,
i.e. \begin{eqnarray}\label{equal}
 \im [H_q(\tilde X, N') \to H_q(\tilde X, N)] = \im
[H'_q(\tilde X) \to H_q(\tilde X, N)].
\end{eqnarray}

For any neighborhood of infinity $N''\subset N'$ one has the following diagram
\begin{eqnarray*}
\begin{array}{ccc}
H_q(\tilde X, N'') & \to & H_q(\tilde X, N')\\ \\
\uparrow & \gamma \searrow& \downarrow  \alpha\\ \\
H'_q(\tilde X) & \underset \beta \to &H_q(\tilde X, N)
\end{array}
\end{eqnarray*}
which gives $\im \beta \subset \im \gamma\subset \im \alpha$ but since we
already know that $\im \alpha$ and $\im \beta$ coincide we obtain $\im \gamma =
\im \alpha$, i.e. (\ref{mittag}).

This completes the proof of Proposition \ref{mittag2} for the specially chosen
neighborhood $N$. If $N_1\subset \tilde X$ is an arbitrary neighborhood of
infinity with respect to $\xi$ then $g_1N\subset N_1$ and we easily see that
for any $N''\subset g_1N'$ one has  $\im \left[H_q(\tilde X,{N}'') \to
H_q(\tilde X,N_1)\right] =\im \left[H_q(\tilde X,g_1N') \to H_q(\tilde
X,N_1)\right]$ i.e. (\ref{mittag}) is satisfied.
\end{proof}

\specialsection*{\bf Part III: Cohomological estimates for $\cat^1(X,\xi)$}

In Part III we combine the results of Parts I and II to obtain new
cohomological lower bounds for $\cat^1(X,\xi)$. This allows us to
compute explicitly $\cat^1(X,\xi)$ is some examples. Finally, we
compare $\cat^1(X,\xi)$ with the values of a similar invariant
$\cat(X,\xi)$ and conclude that their difference can be
arbitrarily large.

\section{Line bundles, algebraic integers and movability of homology classes}

Let $X$ be a finite cell complex and $\xi\in H^1(X;\R)$ be a nonzero cohomology
class. $\xi$ determines the obvious homomorphism $H_1(X;\Z)\to\R$. Its kernel
will be denoted $\ker(\xi)$. The factor-group $H=H_1(X;\Z)/\ker(\xi)$ is a finitely
generated free abelian group which is naturally isomorphic to the group of
periods of $\xi$. The rank of $H$ is equals the rank of class $\xi$; it is
denoted by $r=\rk(\xi)$. Consider the covering $p: \tilde X\to X$ corresponding
to $\ker(\xi)$. This covering has $H$ as the group of covering transformations.

Let $\V_\xi= (\C^\ast)^r= \Hom(H, \C^\ast)$ denote the variety of all complex
flat line bundles $L$ over $X$ such that the induced flat line bundle $p^\ast
L$ over $\tilde X$ is trivial. If $t_1, \dots, t_r\in H$ is a basis, then the
monodromy of $L\in \mathcal V_\xi$ along $t_i$ is a nonzero complex number
$x_i\in \C^\ast$ and the numbers $x_1, \dots, x_r\in \C^\ast$ form a coordinate
system on $\mathcal V_\xi$. Given a flat line bundle $L\in \V_\xi$ the
monodromy representation of $L$ is the ring homomorphism
\begin{eqnarray}\label{monodromy}{\rm
Mon}_L: \Z[H]\to \C\end{eqnarray} sending each $t_i\in H$ to $x_i\in \C^\ast$.

The dual bundle $L^\ast\in \V_\xi$ is such that $L\otimes L^\ast$ is trivial;
if $x_1, \dots, x_r\in \C^\ast$ are coordinates of $L$ then $x_1^{-1}, \dots,
x_r^{-1}\in \C^\ast$ are coordinates of $L^\ast$.

Any nonzero element $P\in \Z[H]$ can be written as $P =\sum_{i=1}^k a_i h_i$
where $a_i\in \Z$, $a_i\not=0$, $h_i\in H$ and $\xi(h_1)<\xi(h_2) < \dots <
\xi(h_k)$. The nonzero integer $a_k$ is called {\it the $\xi$-top coefficient
of $P$.}

The following notion was introduced in \cite{farbe3}, Definition 1.53.

\begin{definition}
A flat line bundle $L\in \V_\xi$ is called a $\xi$-algebraic integer if the
kernel of the monodromy homomorphism ${\rm {Mon}}_L: \Z[H]\to \C^\ast$ contains
a nonzero polynomial $P\in \Z[H]$ having $\xi$-top coefficient 1.
\end{definition}

\begin{theorem}\label{obstruction}
Let $L\in \mathcal V_\xi$ be not a $\xi$-algebraic integer.
 Suppose that for some $v\in H^q(X;L)$ and $z\in H_q(\tilde X;\Z)$ one has $\langle v,
p_\ast(z)\rangle \not=0\in \C$ where $p_\ast: H_q(\tilde X;\Z)\to
H_q(X;L^\ast)$ is the obvious coefficient map. Then the class $z$
is not movable to infinity of $\tilde X$ with respect to $\xi$.
\end{theorem}
\begin{proof} We will show that if a homology class $z\in H_q(\tilde X;\Z)$ is movable to
infinity with respect to $\xi$ then $p_\ast(z)=0\in H_q(X;L^\ast)$ for any
$L\in \V_\xi$ such that the dual bundle $L^\ast$ is not a $\xi$-algebraic
integer. This statement clearly implies the theorem.

Let $S_\xi\subset \Lambda=\Z[H]$ denote the set of all nonzero Laurent
polynomials $P\in \Lambda$ having $\xi$-lowest coefficient 1. The monodromy
homomorphism ${\rm {Mon}}_{L^\ast}: \Lambda\to \C$ is injective when restricted
to $S_\xi$ (because of our assumption that $L$ is not a $\xi$-algebraic
integer). Hence ${\rm {Mon}}_{L^\ast}: \Lambda\to \C$ extends to the
localized ring $\Lambda_\xi=S_\xi^{-1}\Lambda$.

The homomorphism $p_\ast: H_q(\tilde X;\Z)\to H_q(X;L^\ast)$ can be decomposed
as
$$p_\ast: \, H_q(\tilde X;\Z)=H_q(X;\Lambda) \stackrel\alpha\to H_q(X;\Lambda_\xi) \to H_q(X;L^\ast)$$
and the module in the middle equals $H_q(X;\Lambda_\xi)=S_\xi^{-1}H_q(\tilde
X;\Z)$. If $z\in H_q(\tilde X;\Z)$ is movable to infinity with respect to $\xi$
then $\Delta\cdot z=0$ for some $\Delta\in S_\xi$ and hence $\alpha(z)=0$ and
$p_\ast(z)=0$.
\end{proof}

\section{Definition and properties of ${\rm cat}^1(X,\xi)$}

Let $X$ be a finite polyhedron and $\xi\in H^1(X;\R)$ a cohomology class with real coefficients. Let $\omega$ be a closed 1-form on $X$ representing $\xi$, see \cite{farbe2} for the formalism of closed 1-forms on topological spaces.

\begin{definition}
Let $N$ be a positive integer. A subset $A\subset X$ is called $N$-movable with respect to $\omega$, if there exists a continuous homotopy $h_t:A\to X$, $t\in [0,1]$, such that $h_0:A\to X$ is the inclusion and for any point $x\in A$ we have
\begin{eqnarray*}
\int_{h_1(x)}^x \omega & > & N
\end{eqnarray*}
where the integral is calculated along the path $t\mapsto h_{1-t}(x)\in X$, $t\in [0,1]$.
\end{definition}

Recall that for $A\subset X$, ${\rm cat}_X(A)$ denotes the Lusternik-Schnirelmann category of $A$ in $X$, i.e.\ the minimal integer $k$ such that $A$ can be covered by $k$ open sets in $X$ each of which is null-homotopic in $X$.

The following notion has been introduced in \cite{FK}.

\begin{definition}
Let $X$ be a finite polyhedron and $\xi\in H^1(X;\R)$. Fix a closed 1-form $\omega$ in $\xi$. The number ${\rm cat}^1(X,\xi)$ is the minimal integer $k$ such that there exists a closed subset $A\subset X$ with ${\rm cat}_X(X-A)\leq k$ and such that $A$ is $N$-movable with respect to $\omega$ for any positive integer $N$.
\end{definition}

By reversing the order of quantifiers one obtains another notion
originally introduced in \cite{farbe2}.

\begin{definition}\label{defcat}
Let $X$ be a finite polyhedron and $\xi\in H^1(X;\R)$. Fix a closed 1-form $\omega$ in $\xi$. The number ${\rm cat}(X,\xi)$ is the minimal integer $k$ such that for any positive integer $N$ there exists a closed subset $A\subset X$ which is $N$-movable with respect to $\omega$ and such that ${\rm cat}_X(X-A)\leq k$.
\end{definition}

It is easy to see that neither ${\rm cat}^1(X,\xi)$ nor ${\rm
cat}(X,\xi)$ depend on the choice of $\omega$. Furthermore both
notions are homotopy invariants of the pair $(X,\xi)$, see
\cite{farbe2,FK}. Another observation is that for $\xi=0$ we get
the ordinary Lusternik-Schnirelmann category
$\cat(X,\xi)=\cat^1(X,\xi) =\cat(X)$.

It follows straightforwardly from the definitions that
\[ {\rm cat}(X,\xi) \,\, \leq \,\, {\rm cat}^1(X,\xi) \,\, \leq \,\, {\rm cat}(X). \]
We show later in this paper that for some pairs $(X,\xi)$ one has
\begin{eqnarray*}
{\rm cat}(X,\xi) & < & {\rm cat}^1(X,\xi)
\end{eqnarray*}
and that the difference between $\cat^1(X,\xi)$ and $\cat(X,\xi)$
can indeed be arbitrarily large.

\section{The main estimate}

\begin{theorem}\label{main} Let $X$ be a finite cell complex and $\xi\in H^1(X;\R)$.
Let $L\in \V_\xi$ be a complex flat line bundle over $X$ which is
not a $\xi$-algebraic integer. Assume that for some $u\in
H^q(X;L)$ and $z\in H_q(X;\Lambda)$ the evaluation $$\langle u,
p_\ast(z)\rangle \not=0\in \C$$ is nonzero where $p_\ast:
H_q(X;\Lambda) \to H_q(X;L^\ast)$ is the coefficient homomorphism.
Then\footnote{The group $H_q(X;\Lambda)$ is naturally isomorphic
to $H_q(\tilde X;\Z)$. However the category weights of $z$ viewed
as element of $H_q(X;\Lambda)$ or of $H_q(\tilde X;\Z)$ are in
general different. In inequality (\ref{mainestim}) the symbol
$\cwgt(z)$ denotes the category weight of $z$ regarded as an
element of $H_q(X;\Lambda)$.}
\begin{eqnarray}\label{mainestim}
\cat^1(X,\xi)\geq \cwgt(z) +1.
\end{eqnarray}
\end{theorem}
\begin{proof}
Denote $k=\cwgt(z)$ and assume the contrary, i.e.that $\cat^1(X,\xi)\leq k$.
Then there exists a closed subset $A\subset X$ with $\cat_XA\leq k$ such that
the complement $F=X-A$ is $N$-movable for any $N>0$ with respect to a closed
1-form $\omega$ on $X$ representing $\xi$. Applying the definition, we find
that $z$ can be realized by a singular cycle $c$ in $X-A=F$ with coefficients
in the local system $\Lambda$.

Consider the covering $p: \tilde X\to X$ corresponding to $\ker(\xi)$. Viewed
differently, the cycle $c$ is a usual singular cycle in $\tilde X$ lying in the
set $\tilde F=p^{-1}(F)$. Since $F$ is $N$-movable for any $N$ we find that any
cycle in $\tilde F$ is movable to infinity with respect to $\xi$. Thus we
obtain a contradiction with Theorem \ref{obstruction}.
\end{proof}

\begin{theorem}\label{main2}
Let $X$ be a finite cell complex and $\xi\in H^1(X;\R)$. Let $L\in
\V_\xi$ be a complex flat line bundle over $X$ which is not a
$\xi$-algebraic integer. Suppose that for an integral homology
class $z\in H_q(\tilde X;\Z)=H_q(X;\Lambda)$ and some cohomology
classes $u\in H^{d}(X;L)$ and $u_i\in H^{d_i}(X;\C)$, where
$d_i>0$ for $i=1, \dots, k$, the evaluation $\langle u\cup
u_1\cup\dots\cup u_k, p_\ast(z)\rangle \not=0\in \C$ is nonzero.
Here $p_\ast(z)\in H_d(X;L^\ast)$, $q=d+d_1+\dots+d_k$. Then one
has
\begin{eqnarray}
\cat^1(X,\xi)\geq \cwgt(z)+k +1. \end{eqnarray}
\end{theorem}
Here $\cwgt(z)$ denotes the category weight of $z$ viewed as a homology class
of $X$ with local coefficient system $\Lambda$.
\begin{proof} First observe that we may assume that the classes $u_1, \dots,
u_k$ are integral, i.e. lie in $H^\ast(X;\Z)$. Indeed, the product
$\langle u\cup u_1\cup\dots\cup u_k, p_\ast\rangle$ is a
multilinear function of $u_1, \dots, u_k$; since the integral
cohomology classes generate $H^\ast(X;\C)$ vanishing of this
function on all integral combinations would imply vanishing in
general.

 Denote $z'=p^\ast(u_1\cup\dots\cup u_k)\cap z\in H_d(\tilde X;\Z)=H_q(X;\Lambda)$. Then
 $$\langle u, p_\ast(z')\rangle = \langle u\cup u_1\cup\dots\cup u_k,
p_\ast(z)\rangle\not=0\in \C.$$ Applying the previous theorem we
find $\cat^1(X,\xi)\geq \cwgt(z')+1$. Now, Theorem \ref{cap1}
gives $\cwgt(z')\geq k+\cwgt(z)$. This completes the proof.
\end{proof}

\noindent{\bf Remark:} Consider the statement of Theorem
\ref{main2} in the special case $\xi=0$. Then the variety $\V_\xi$
contains the trivial line bundle $L=\C$ only and $L=\C$ is not a
$\xi$-algebraic integer. Hence Theorem \ref{main2} gives the
inequality
$$\cat(X)\geq \cwgt(z)+k+1$$ under the assumption that
$$\langle u_1\cup\dots\cup u_k,z\rangle \not=0$$
where $u_i\in H^{d_i}(X;\C)$, $d_i>0$ and $z\in H_d(X;\C)$,
$d=d_1+\dots+d_k$. This claim is a special case of
(\ref{improved}).

\begin{example}\label{firstinst}
 Let $X=\Sigma$ be a closed orientable surface of genus
$g>1$ and $\xi\not=0\in H^1(X;\R)$. Fix a flat line bundle $L\in \V_\xi$ which
is transcendental, see \cite{FS1}, \S 6. Then $H^1(X;L)$ has dimension
$2g-2>0$. Pick a nonzero class $u\in H^1(X;L)$. By Proposition 6.5 from
\cite{FS1} there exists a homology class $z\in H_1(X;\Lambda)$ such that
$\langle u, p_\ast(z)\rangle \not=0$. Since $\cwgt(z) \geq \swgt(z)\geq 1$ we
get $\cat^1(\Sigma, \xi)\geq 2$ by applying Theorem \ref{main}. Since
$\cat^1(X,\xi)\leq \dim X$ in general for $\xi\not=0$ we find
\begin{eqnarray}
\cat^1(\Sigma, \xi)=2
\end{eqnarray}
for any nonzero $\xi\in H^1(\Sigma;\R)$.

Note that $\cat(\Sigma, \xi)=1$ for any $\xi\not=0$, see Theorem 12 in
\cite{FS1}. This gives a first instance where
\begin{eqnarray}
\cat(X,\xi) < \cat^1(X,\xi).
\end{eqnarray}
\end{example}

\section{A controlled version of ${\rm cat}^1(X,\xi)$}

We have seen in Example \ref{firstinst} that ${\rm cat}(X,\xi)$
and ${\rm cat}^1(X,\xi)$ can indeed be different. In order to show
that the difference between them can be arbitrary large, we have
to introduce a controlled version of ${\rm cat}^1(X,\xi)$ which
behaves better under cartesian products. The following discussion
is very similar to \cite[Section 9]{FS1}.

Let $\omega$ be a continuous closed 1-form on a finite cell complex $X$. Let
$\xi=[\omega]\in H^1(X;\R)$ be the cohomology class represented by $\omega$.

\begin{definition}\label{cmovable}
Let $N$ and $C$ be two positive integers. A subset $A\subset X$ is $N$-movable with
respect to $\omega$ with control $C$ if there exists a continuous homotopy
$h_t:A\to X,$ $t\in [0,1],$ such that (1) $h_0:A\to X$ is the inclusion; (2)
for any point $x\in A$ one has
\begin{eqnarray}
\int\limits_{x}^{h_1(x)} \omega\,<\,- N,\label{cint}
\end{eqnarray}
where the integral is calculated along the path $t\mapsto h_{t}(x)\in X$, $t\in
[0,1]$ and (3) for any point $x\in A$ and for any $t\in [0,1]$ one has
\begin{eqnarray}
\int\limits_{x}^{h_t(x)} \omega\,\leq\,C.\label{ccint}
\end{eqnarray}
\end{definition}

\begin{definition}\label{cdeff3}
Fix a closed 1-form $\omega$ representing $\xi$. The number $\ccat^1(X,\xi)$ is the minimal integer $k$ with the property that there exists $C>0$ and a closed subset $A\subset X$ with $\cat_X(X-A)\leq k$ and such that $A$ is $N$-movable with control $C$ with respect to $\omega$ for every positive integer $N$.
\end{definition}

\begin{lemma}\label{lm16}
The following properties hold for ${\rm ccat}^1(X,\xi)$.
\begin{enumerate}
\item We have ${\rm cat}^1(X,\xi)\leq {\rm ccat}^1(X,\xi)$.
\item If $X$ is connected and $\xi\not=0$, then ${\rm ccat}^1(X,\xi)\leq {\rm cat}(X)-1$.
\item If $\xi=0$, then ${\rm ccat}^1(X,\xi)={\rm cat}(X)$.
\item If $\phi:Y\to X$ is a homotopy equivalence and $\xi\in H^1(X;\R)$, then
\begin{eqnarray*}
{\rm ccat}^1(X,\xi)&=&{\rm ccat}^1(Y,\phi^\ast \xi)
\end{eqnarray*}
\end{enumerate}
\end{lemma}

\begin{proof}
The first assertion is obvious, the remaining assertions are obtained by repeating the arguments given in \cite{farbe2} and \cite{FS1}.
\end{proof}

\begin{remark}
It is worth pointing out that the applications of ${\rm cat}^1(X,\xi)$ to dynamics described in \cite{FK} also hold with the potentially larger quantity ${\rm ccat}^1(X,\xi)$, compare \cite[Remark 9.9]{FS1}.
\end{remark}

The desired product inequality now reads as follows.

\begin{theorem}\label{prodineq}
Let $X$ and $Y$ be finite cell complexes and let $\xi_X\in H^1(X;\R)$ and
 $\xi_Y\in H^1(Y;\R)$ be real cohomology classes. Assume that
\begin{eqnarray}\label{posi}
\ccat^1(X, \xi_X)>0\quad \mbox{or}\quad \ccat^1(Y, \xi_Y)>0.
\end{eqnarray}
 Then
\begin{eqnarray}\label{ineqprod}
\ccat^1(X\times Y, \xi) \leq \ccat^1(X, \xi_X) + \ccat^1(Y, \xi_Y) -1,
\end{eqnarray}
where \begin{eqnarray}\xi=\xi_X\times 1 \, +\,  1\times \xi_Y.\end{eqnarray}
\end{theorem}

We skip the proof since it is fully analogous to the proof of the
similar statement for $\ccat(X,\xi)$ given in \cite[Theorem
9]{FS1}.

\section{Calculation of $\cat^1(X,\xi)$ for products of surfaces}

\begin{theorem}\label{surface}
Let $M^{2k}$ denote the product $\Sigma_1\times \Sigma_2\times
\dots\times \Sigma_k$ where each $\Sigma_i$ is a closed orientable
surface of genus $g_i>1$. Given a cohomology class $\xi\in
H^1(M^{2k};\R)$, one has
\begin{eqnarray}
\cat^1(M^{2k},\xi)=\ccat^1(M^{2k},\xi) = 1+k+ r \end{eqnarray} where $r$
denotes the number of indices $i\in \{1, 2, \dots, k\}$ such that the
cohomology class $\xi|_{\Sigma_i}\in H^1(\Sigma_i,\R)$ vanishes. In particular
\begin{eqnarray}\label{one2}
\cat^1(M^{2k},\xi)=\ccat^1(M^{2k},\xi) = 1+k
\end{eqnarray}
assuming that $\xi|_{\Sigma_i}\not=0\in H^1(\Sigma_i;\R)$ for any $i=1, \dots,
k$.
\end{theorem}

\begin{proof}
After rearranging the surfaces we may assume that
$\xi_i=\xi|_{\Sigma_i}$ is nonzero for $i=1,\ldots, k-r$ and
$\xi_i=0$ for $i>k-r$.

Note that $\ccat^1(\Sigma_i,\xi_i)>0$ for any $i=1, \dots, k$.
Indeed, otherwise applying Theorem 10 of \cite{FS1} we would get
$\chi(\Sigma_i)=0$ contradicting our assumption $g_i>0$. Hence we
may apply the inequality of Theorem \ref{prodineq} several times
to obtain
\begin{eqnarray*}
\ccat^1(M^{2k},\xi)&\leq & \sum_{i=1}^k \ccat^1(\Sigma_i,\xi_i)
-(k-1).
\end{eqnarray*}
By Example \ref{firstinst} and Lemma \ref{lm16} we have
\begin{eqnarray*}
\ccat^1(\Sigma_i,\xi_i)= \cat^1(\Sigma_i,\xi_i)  =
\left\{\begin{array}{ccl} 2&{\rm
if}&i\leq k-r,\\ \\
3&{\rm if}&i>k-r,\end{array} \right.
\end{eqnarray*}
and thus
\begin{eqnarray}\label{upperbound}
\ccat^1(M^{2k},\xi)&\leq & 2(k-r)+3r-(k-1)\, =\,k+r+1.
\end{eqnarray}

Next we prove the opposite inequality to (\ref{upperbound}). Let
$L\in \V_\xi$ be transcendental. Denote $H=\pi_1(M)/\ker\, \xi$
and $L_i=L|_{\Sigma_i}$ and $H_i=\pi_1(\Sigma_i)/\ker \,\xi_i$. It
follows that $L_i$ is also transcendental. Choose $u_i'\in
H^1(\Sigma_i;L_i)$ and $z_i\in H_1(\tilde{\Sigma}_i;\Z)$ such that
$\langle u_i',p_\ast(z_i)\rangle \not=0$ as in Example
\ref{firstinst}. Here $p_i: \tilde{\Sigma}_i\to \Sigma_i$ is the
covering space corresponding to $\ker\,\xi_i$ and ${p_i}_\ast:
H_\ast(\tilde \Sigma_i,\Z) \to H_i(\Sigma_i, L_i^\ast)$. Note that
for $i>k-r$ we simply have $\tilde{\Sigma}_i=\Sigma_i$ and
$L_i=\C$. Now, $\ker\,\xi_1\times \ldots\times \ker\,\xi_k\subset
\ker \,\xi$ so there is a covering map
$$q:\tilde{\Sigma}_1\times \ldots \times \tilde{\Sigma}_k\to
\tilde{M}$$ where $\tilde{M}$ is the covering space of $M$
corresponding to $\ker \,\xi$. Let
\[
z'=z_1\times \ldots \times z_k \in H_k(M;\Z[H_1\times \ldots \times H_k])
\cong H_k(\tilde{\Sigma}_1\times \ldots \times \tilde{\Sigma}_k;\Z)
\]
and
\begin{eqnarray*}
z&=&q_\ast(z')\,\,\in \,\, H_k(M;\Z[H])\,\,\cong \,\, H_k(\tilde{M};\Z).
\end{eqnarray*}
It follows from Corollary \ref{prodweight} and Lemma
\ref{coeffchange} that $\cwgt(z)\geq k$ (where $z$ is viewed as an
element of $H_k(M;\Z[H])$).

Denote
\begin{eqnarray*}
u &=&u'_1\times \ldots \times u'_{k-r} \times 1 \times \ldots
\times 1 \,\,\in \,\, H^{k-r}(M;L)
\end{eqnarray*}
and
\begin{eqnarray*}
u_j&=&p^\ast_{k-r+j}u'_{k-r+j}\,\,\in\,\,H^1(M;\C), \quad
j=1,\ldots,r,
\end{eqnarray*}
where $p_{k-r+j}:M\to \Sigma_{k-r+j}$ is projection. Notice that
\begin{eqnarray*}
\langle u\cup u_1\cup\ldots\cup u_r,p_\ast(z)\rangle &=& \pm \,
\prod_{i=1}^k \langle u'_i,{p_i}_\ast(z_i)\rangle \,\,\not= \,\,0.
\end{eqnarray*}
Theorem \ref{main2} and Corollary \ref{prodweight} apply and give
\begin{eqnarray*}
\cat^1(M,\xi) &\geq & \cwgt (z) +r +1\\ \\
& \geq & k+r+1.
\end{eqnarray*}
Combining this with (\ref{upperbound}) we obtain
\[
\cat^1(M,\xi)\,\,=\,\,\ccat^1(M,\xi)\,\,=\,\, k+r+1 \]
as claimed.
\end{proof}

We now want to compare the values of $\cat^1(M,\xi)$ with the
invariant $\cat(M,\xi)$ (see Definition \ref{defcat}) for products
of surfaces $M= \Sigma_1\times \dots\times \Sigma_k$ where each
$\Sigma_i$ is a closed orientable surface of genus $g_i>1$. It was
shown in \cite[Thm.17]{FS1} that one has
\begin{eqnarray}
\cat(M,\xi)&=&1+2r
\end{eqnarray}
where $r$ denotes the number of indices $i\in \{1, \dots, k\}$
such that $\xi|_{\Sigma_i}=0$.

\begin{corollary}\label{cor} Under the assumptions of Theorem \ref{surface} the difference
\begin{eqnarray}\label{difference}
\cat^1(M,\xi)-\cat(M,\xi)\end{eqnarray} equals the number of
indices $i\in \{1, \dots, k\}$ such that $\xi|_{\Sigma_i}\not=0\in
H^1(\Sigma_i;\R)$.
\end{corollary}

Corollary \ref{cor} leads to the following statement which is one
of the main results of this paper:

\begin{corollary}
The difference (\ref{difference}) can be arbitrarily large.
\end{corollary}

\bibliographystyle{amsalpha}

\end{document}